\documentclass[12pt,twoside,a4paper]{amsart}
\usepackage{amssymb}
\usepackage{amsthm}
\usepackage[english]{babel}
\usepackage[latin1]{inputenc} 
\usepackage[T1]{fontenc}

\def\End{{\rm End}}
\def\deg{\mathrm{deg}\,}

\def\z{\zeta}

\def\deta{\delta_{\eta}}

\def\w{\wedge}
\def\db{\bar\partial}
\def\dbar{\bar\partial}
\def\ddbar{\partial\dbar}

\def\C{{\mathbb C}}
\def\P{{\mathbb P}}

\def\Cn{\C^n}
\def\Pn{\P^n}
\def\Pm{\P^m}

\def\Re{{\rm Re\,  }}

\def\End{\textrm{End}}
\def\Span{\textrm{Span}}
\def\be{\begin{equation}}
\def\ee{\end{equation}}

\newtheorem{thm}{Theorem}[section]
\newtheorem{lma}[thm]{Lemma}

\newtheorem{prop}[thm]{Proposition}

\newtheorem{De}{Definition}
\newcommand{\df}[1] {\begin{De}\emph{#1}\end{De}}

\newtheorem{preremark}{Remark}
\newtheorem{preex}{Example}

\begin{document}

\title{Weighted integral formulas on manifolds}
\author{Elin Götmark}

\maketitle

\begin{abstract}
We present a method of finding weighted Koppelman formulas for
$(p,q)$-forms on $n$-dimensional complex manifolds $X$ which admit a
vector bundle of rank $n$ over $X \times X$,
such that the diagonal of $X \times X$ has a defining section. We
apply the method to $\Pn$ and find weighted Koppelman formulas for
$(p,q)$-forms 
with values in a line bundle over $\Pn$. As an application, we
look at the cohomology groups of $(p,q)$-forms over $\Pn$ with values
in various line bundles, and find explicit solutions to the
$\dbar$-equation in some of the trivial groups. We also look at
cohomology groups of $(0,q)$-forms over $\Pn \times \Pm$ with values in
various line bundles. Finally, we apply our method to developing weighted
Koppelman formulas on Stein manifolds.
\end{abstract}

\tableofcontents

\section{Introduction}

The Cauchy
integral formula provides a decomposition of a holomorphic function in
one complex variable in simple rational functions, and is a
cornerstone in function theory in one complex variable. The kernel is
holomorphic and works for any domain. In several complex variables it
is harder to find appropriate representations. The simplest
multivariable analog, the Bochner-Martinelli kernel, is not as useful
since the kernel is not holomorphic. The Cauchy-Fantappie-Leray formula
is a generalization which gives a holomorphic kernel in domains which
admit a holomorphic support function. Henkin and Ramirez in \cite{HE},
\cite{RA} obtained 
holomorphic kernels in strictly pseudoconvex domains $G$ by finding such
support functions. Henkin also found solutions to the $\dbar$-equation in
such domains. This was done by means of a Koppelman formula, which
represents a $(p,q)$-form $\phi$ defined in some domain $D$ as a sum
of integrals  
\[
\phi(z) =  \int_{\partial D} K \w \phi + \int_{D} K \w \dbar \phi +  
\dbar_z \int_{D} K \w \phi + \int_{D} P \w \phi,
\]
by means of the current $K$ and the smooth form $P$. If $\phi$ is a closed 
form and the first and fourth terms of the right hand side 
of Koppelman's formula vanish, we get a 
solution of the $\dbar$-problem for $\phi$.
Henkin's result paved the way for the Henkin-Skoda theorem
(see \cite{HE2} and \cite{SK}), which provided improved
$L^1$-estimates on $\partial G$ for solutions of the $\dbar$-equation
by weighting the integral formulas.  

Andersson and Berndtsson \cite{AB} found a flexible method of
generating weighted formulas for representing holomorphic functions
and solutions of the $\dbar$-equation. It was further developed by
Berndtsson \cite{BE3} to find solutions to division and interpolation
problems. If $V$ is a regular analytic subvariety of some domain $D$
in $\Cn$ and $h$ is holomorphic in $V$, then Berndtsson found a kernel
$K$ such that 
\[
H(z) = \int_V h(\z) K(\z,z)
\]
is a holomorphic function which extends $h$ to $D$. 
If $f = (f_1, \ldots, f_m)$ are holomorphic functions without 
common zeros, he also found a solution to the division
problem $\phi = f \cdot p$ for a given
holomorphic function $\phi$. Passare \cite{PA} used
weighted integral formulas to solve a similar division problem, where
the $f_i$'s do have common zeros, but the zero sets have a complete
intersection. He also proved the duality theorem for complete
intersections (also proved independently by Dickenstein and Sessa
\cite{DS}). Since then weighted
integral formulas have been used by a row of authors to obtain
qualitative estimates of solutions of the $\dbar$-equation and of
division and interpolation problems, for example sharp
approximation by polynomials \cite{ZE},
estimates of solutions to the Bézout equation \cite{BY}, and explicit
versions of the fundamental principle \cite{BP}. More
examples and references can be found in the book \cite{BY2}. More recently, 
Andersson \cite{MA3} introduced a method generalizing \cite{AB} and 
\cite{BE3} which is even more flexible and also easier to handle. It
allows for some recently found representations with residue currents,
applications to division and interpolation problems, and also allows
for $f$ to be a matrix of functions. 

There have been several attempts to obtain integral formulas on
manifolds. Berndtsson \cite{BE1} gave a method of obtaining integral
kernels on $n$-dimensional manifolds $X$ which admit a vector bundle
of rank $n$ over $X \times X$
such that the diagonal has a defining section, but did not consider
weighted formulas. Formulas on Stein manifolds were treated first in
Henkin and Leiterer \cite{HL}, where formulas for $(0,q)$-forms are
found, then in Demailly and Laurent-Thiébaut \cite{DLT}, where the
leading term in a kernel for $(p,q)$-forms is found, in Andersson
\cite{MA4}, which is a generalization of \cite{AB} following Henkin
and Leiterer, and finally in Berndtsson \cite{BE1} where the
method described therein is applied to Stein manifolds. Formulas on
$\Pn$ have been considered in \cite{PH}, where they 
were constructed by using known formulas in $\mathbb{C}^{n+1}$, and in
\cite{BE2}, where they were constructed directly on $\Pn$. There is
also an example at the end of Berndtsson \cite{BE1}, where the
method of that article is applied to $\Pn$. 

In this article, we begin in Section \ref{regn} by developing a method
for generating weighted  
integral formulas on $\Cn$, following \cite{MA1}. Section \ref{yuuri}
describes a similar method which can be used on $n$-dimensional
manifolds $X$ which admit a vector bundle  
of rank $n$ over $X \times X$ such that the diagonal has a defining
section. It has similar results as the method described 
in \cite{BE1}, but with the added benefit of yielding weighted
formulas. The method of 
Section \ref{yuuri} is applied to complex projective space $\Pn$
in Section \ref{svamp}, where we find a Koppelman formula for
differential forms with values in a line bundle over $\Pn$. In the
$\Pn$ case we get formulas which coincide with Berndtsson's formulas
in \cite{BE2} in the case $p=0$, but they are not the same in the general
$(p,q)$-case.

As an application, in Section \ref{wombat} we
look at the cohomology groups of $(p,q)$-forms over $\Pn$ with values
in various line bundles, and find which of them are trivial (though we
do not find all the trivial groups). Berndtsson's formulas
in \cite{BE2} give the same result. The trivial cohomology groups of
the line bundles over $\Pn$ are, of course,
known before, but our method gives explicit solutions of the
$\dbar$-equations. In Section \ref{sol} we look instead at cohomology
groups of $(0,q)$-forms over $\Pn \times \Pm$ with values in various
line bundles. Finally, in Section \ref{sjunga} we apply the method of
Section \ref{yuuri} to finding weighted integral formulas on Stein
manifolds, following \cite{HL} but also developing weighted formulas.

\section{Weighted Koppelman formulas in $\Cn$}
\label{regn}

As a model for obtaining representations on manifolds, we present the
$\Cn$ case in some detail. The material in this section follows the
last section of 
\cite{MA1}. The article \cite{MA1} is mostly concerned with
representation of holomorphic functions, but in the last 
section a method of constructing weighted Koppelman formulas in
$\Cn$ is indicated. We expand this material and give proofs in more
detail. We begin with some motivation from the one-dimensional case:

One way of obtaining a representation formula for a holomorphic
function would be to solve the equation  
\[
\dbar u = [z],
\]
where $[z]$ is the Dirac measure at $z$ considered as a
$(1,1)$-current, since then one would get 
an integral formula by Stokes' theorem. Less 
obviously, note that the kernel of Cauchy's integral formula in
$\mathbb{C}$ also satisfies the equation 
\[
\delta_{\z-z} u = 1,
\]
where $\delta_{\z-z}$ denotes contraction with the vector field $2 \pi i
(\z-z) \partial /\partial \z$. These two can be combined into the
equation 
\be \label{zira}
\nabla_{\z-z} u := (\delta_{\z-z} - \dbar) u = 1-[z].
\ee
To find representation formulas for holomorphic functions in $\Cn$, we
look for solutions to equation (\ref{zira}) in $\Cn$, where
$\delta_{\z-z}$ is contraction with 
\[
2 \pi i\sum (\z_j-z_j) \frac{\partial}{\partial \z_j}.
\]
Since the right hand side of (\ref{zira}) contains one form of
bidegree $(0,0)$ and one of bidegree $(n,n)$, we must in fact have $u
= u_{1,0} + u_{2,1} + \ldots + u_{n,n-1}$, where $u_{k,k-1}$ has
bidegree $(k,k-1)$. We can then write (\ref{zira}) as the system 
of equations 
\[
\delta_{\z-z} u_{1,0} = 1, \qquad \delta_{\z-z} u_{1,2} - \dbar u_{1,0} = 0 \qquad 
\ldots \qquad \dbar u_{n,n-1} = [z]. 
\]
In that case, $u_{n,n-1}$ will
satisfy $\dbar u_{n,n-1} = [z]$ and will give a kernel for a
representation formula.  
The advantage of this approach is that it easily 
allows for weighted integral formulas, as we will see. 

To get Koppelman formulas for $(p,q)$-forms, we need
to consider $z$ as a variable and not a constant. If we find
$u_{n,n-1}$ such that $\dbar u_{n,n-1} = [\Delta]$, where $\Delta =
\{(\z,z): \z = z\}$ is
the diagonal of $\Cn_\z \times \Cn_z$ and $[\Delta]$ is the current of
integration over $\Delta$, then $u_{n,n-1}$ will be the kernel that we
seek. In fact, if we let $\phi$
be a $(p,q)$-form, and $\psi$ an $(n-p,n-q)$ test form, we have 
\[
\int_z \left( \int_\z \phi(\z) \w [\Delta] \right) \w \psi(z) = 
\int_{z,\z} \phi(\z) \w \psi(z) \w [\Delta] = \int_z \phi(z) \w \psi(z)
\]
so that $\int_\z \phi(\z) \w [\Delta] = \phi(z)$ in the current
sense. 


In more detail, then: Let $\Omega$ be a domain in $\Cn$
and let $\eta(\z,z) = 2 \pi i(z-\z)$, where $(\z,z) \in \Omega \times
\Omega$. Note that $\eta$ vanishes to the first order on the diagonal.  
Consider the subbundle $E^\ast = \textrm{Span}\{d\eta_1, \ldots , d
\eta_n \}$ of the cotangent bundle $T^\ast_{1,0}$
over $\Omega \times \Omega$. Let $E$ be its dual bundle, and let 
$\delta_\eta$ be an operation on $E^\ast$, defined as contraction with
the section  

\be \label{toste}
\sum_1^n \eta_j e_j,
\ee
where $\{e_j\}$ is the dual basis to $\{d \eta_j\}$. Note that $\deta$
anticommutes with $\dbar$. 

Consider the bundle $\Lambda (T^\ast(\Omega \times \Omega) \oplus
E^\ast)$ over $\Omega \times \Omega$. An example of an element of the
fiber of this bundle at $(\z,z)$ is $d\z_1 \w d \bar z_2 \w d \eta_3$. 
We define 
\be \label{fisk}
\mathcal{L}^m = \bigoplus_p C^\infty (\Omega \times \Omega,
\Lambda^p E^\ast \w \Lambda^{p+m} T_{0,1}^\ast(\Omega \times \Omega)).
\ee
Note that $\mathcal{L}^m$ is a subset of
the space of sections of $\Lambda (T^\ast(\Omega \times \Omega) \oplus
E^\ast)$. Let $\mathcal{L}_{curr}^m$ be the
corresponding space of currents. If $f \in \mathcal{L}^m$ and $g \in
\mathcal{L}^{k}$, then $f \w g \in \mathcal{L}^{m+k}$. 

We define the operator 

\[
\nabla = \nabla_{\eta} = \delta_\eta - \dbar,
\]
which maps $\mathcal{L}^m$ to $\mathcal{L}^{m+1}$. We see that $\nabla$ 
obeys Leibniz' rule, that is, 
\be \label{medlem}
\nabla (f \w g) = \nabla f \w g + (-1)^m f \w \nabla g,
\ee
if $f \in \mathcal{L}^m$. Note that $\nabla^2 = 0$, which means that 
\[
\ldots \stackrel{\nabla}{\to} \mathcal{L}^m \stackrel{\nabla}{\to}
\mathcal{L}^{m+1} \stackrel{\nabla}{\to} \ldots 
\]
is a complex. We also have the following useful property: If 
$f$ is a form of bidegree $(n,n-1)$ and $D \subset \Omega 
\times \Omega$, then 

\be \label{gröt}
\int_{\partial D} f = - \int_D \nabla f.
\ee
This follows from Stokes' theorem and the fact that $\int_D
\deta f=0$. The operator $\nabla$ is defined also for currents, since
$\dbar$ is defined for currents, and $\deta$ just amounts to 
multiplying with a smooth function, which is also defined for a
current. 

As in the beginning of this section, we want to find a
solution to the equation 
\be \label{bombay} 
\nabla_{\eta} u = 1-[\Delta].
\ee
with $u \in \mathcal{L}_{curr}^{-1}$ (since the left hand side lies in 
$\mathcal{L}_{curr}^{0}$), so as before, we have $u = u_{1,0} +
u_{2,1} + \ldots + u_{n,n-1}$, where $u_{k,k-1}$ has degree $k$ in $E^\ast$
and degree $k-1$ in $T_{0,1}^\ast$. 
\begin{prop}
Let 
\[
b(\z,z) = \frac{1}{2 \pi i} \frac{\partial |\eta|^2}{|\eta|^2}
\]
and 
\be \label{fondjävel}
u_{BM} = \frac{b}{\nabla_\eta b} = \frac{b}{1-\dbar b} = b + b \w \dbar b +
\ldots b \w (\dbar b)^{n-1},
\ee
where we get the right hand side by expanding the fraction in a
geometric series. Then $u$ solves equation (\ref{bombay}).
\end{prop}
\noindent
The crucial step in the proof is showing that $\dbar (b \w (\dbar
b)^{n-1}) = [\Delta]$, which is common knowledge, since $b \w (\dbar
b)^{n-1}$ is actually the well-known Bochner-Martinelli kernel. 

A form $u$ which satifies $\nabla_\eta
u = 1$ outside $\Delta$ is a good candidate for solving equation
(\ref{bombay}). The following proposition gives us a
criterion for when such a $u$ in fact is a solution:
\begin{prop} \label{isis}
Suppose $u \in \mathcal{L}^{-1}(\Omega \times \Omega \setminus 
\Delta)$ solves $\nabla_\eta
u = 1$, and that $|u_k| \lesssim |\eta|^{-(2k-1)}$. We then have  
$\nabla_{\eta} u = 1-[\Delta]$. 
\end{prop}

\begin{proof}
Let $u_{BM}$ be the form defined by (\ref{fondjävel}), and let $u$ 
be a form satisfying the conditions in the proposition. We know that 
$\nabla (u \w u_{BM}) = u_{BM} - u$ pointwise outside $\Delta$, in
light of (\ref{medlem}). We want to show that this also holds in the
current sense, i. e. 

\be 
\int \nabla (u \w u_{BM}) \w \phi = \int (u_{BM} - u) \w \phi,
\ee
where $\phi$ is a test form in $\Omega \times \Omega$. Using firstly
that $u \w u_{BM}$ is locally integrable (since $u \w u_{BM} 
= \mathcal{O}(|\eta|^{-(2n-2)})$ near $\Delta$), and secondly
(\ref{gröt}), we get 
\begin{eqnarray} \label{mikroskop}
& & \int \nabla (u \w u_{BM}) \w \phi = \lim_{\epsilon \to 0}
\int_{|\eta|>\epsilon} (u \w u_{BM}) \w  
\nabla \phi = \nonumber \\
& = & - \lim_{\epsilon \to 0} \left( \int_{|\eta|=\epsilon} u \w u_{BM} 
\w \phi + \int_{|\eta|>\epsilon} \nabla (u \w u_{BM}) \w \phi \right) .
\end{eqnarray}
The boundary integral in (\ref{mikroskop}) will converge to zero 
when $\epsilon \to 0$, since $u \w u_{BM} = 
\mathcal{O}(|\eta|^{-2n+2})$ and $\textrm{Vol}(\{|\eta|=\epsilon \}
\cap \textrm{supp} (\phi)) = 
\mathcal{O}(\epsilon^{2n-1})$. As for the last integral in
(\ref{mikroskop}), we get 

\[
\lim_{\epsilon \to 0} \int_{|\eta|>\epsilon} \nabla (u \w u_{BM}) \w
\phi = \lim_{\epsilon \to 0} \int_{|\eta|>\epsilon} 
(u - u_{BM}) \w \phi = \int (u - u_{BM}) \w \phi,
\]
since $u - u_{BM}$ is locally integrable, thus $\nabla (u \w u_{BM}) =
u_{BM} - u$ as currents.  
It follows that $\nabla u = \nabla u_{BM}$ since $\nabla^2 = 0$, and
since $u_{BM}$ satisfies the equation (\ref{bombay}), $u$ must also do so. 
\end{proof}

\begin{preex} \label{wayward}
\emph{If $s$ is a smooth $(1,0)$-form in $\Omega \times \Omega$ such that 
$|s| \lesssim |\eta|$ and $|\deta s| \gtrsim |\eta|^2$, we can set 
$u = s/\nabla s$. By Proposition \ref{isis}, $u$ will satisfy equation
(\ref{bombay}), and 
\[
u_{n,n-1} = \frac{s \w (\dbar s)^{n-1}}{(\deta s)^n}
\]
is the classical Cauchy-Fantappie-Leray kernel.}
\end{preex}

We now introduce weights, which will allow us to get more flexible
integral formulas:

\df{A form $g \in \mathcal{L}^{0}(\Omega \times \Omega)$ is a 
weight if $g_{0,0}(z,z) = 1$ and $\nabla_{\eta}g = 0$.} 
\noindent
The form $1 + \nabla Q$ is an example of a weight, if 
$Q \in \mathcal{L}^{-1}$. In fact, we have considerable flexibility
when choosing weights: if $g = 1 + \nabla Q$ is a weight and $G(\lambda)$ is a
holomorphic function such that $G(0) = 1$, then it is easy to see that 
\[
G(g) = \sum_0^n G^{(k)}(\deta Q)(-\dbar Q)^k/k!
\]
is also a weight. We can now prove the following representation formula:

\begin{thm} [Koppelman's formula] \label{koala}
Assume that $D \subset \subset \Omega$, $\phi \in \mathcal{E}_{p,q}
(\bar D)$, and that the current $K$ and the smooth form $P$ solve the
equation  
\be \label{miyazaki}
\dbar K = [\Delta] - P.
\ee
We then have 
\be \label{iris}
\phi(z) =  \int_{\partial D} K \w \phi + \int_{D} K \w \dbar \phi +  
\dbar_z \int_{D} K \w \phi + \int_{D} P \w \phi,
\ee
where the integrals are taken over the $\z$ variable. 
\end{thm}
\begin{proof}
First assume that $\phi$ has compact support in $D$, so that the 
first integral in (\ref{iris}) vanishes. Take a test form $\psi(z)$ of
bidegree $(n-p,n-q)$ in $\Omega$. Then we have

\begin{eqnarray}
& & \int_z \left( \int_\z K \w \dbar \phi + 
\dbar_z \int_\z K \w \phi + \int_{\z} P \w \phi \right) \w \psi = 
\nonumber \\
& = & \int_{z,\z} K \w d \phi \w \psi + 
(-1)^{p+q} \int_{z,\z} K \w \phi \w d \psi + \int_{z,\z} P \w \phi \w \psi = 
\nonumber \\
& = & \int_{z,\z} K \w d (\phi \w \psi) + \int_{z,\z} P \w \phi 
\w \psi = \nonumber \\
& = & \int_{z,\z} dK \w \phi \w \psi + \int_{z,\z} P \w 
\phi \w \psi = \int_z \phi \w \psi, \nonumber 
\end{eqnarray}
where we use Stokes' theorem repeatedly. If $\phi$ does not have compact 
support in $D$, we can prove the general case e g by making the
decomposition $\phi = \phi_1 +  
\phi_2$, where $\phi_1$ has compact support in $D$, and $\phi_2(\z) = 0$ 
in a neighborhood of $z$. 
\end{proof}
It is easy to obtain $K$ and $P$
which solve (\ref{miyazaki}): If we take $g$ to be a weight and $u$ to
be a solution of (\ref{bombay}), then we can solve the equation 
\[
\nabla_{\eta} v = g - [\Delta]
\]
by choosing $v = u \w g$. This means that $K = (u \w g)_{n,n-1}$ and $P
= g_{n,n}$ will solve (\ref{miyazaki}). 

\begin{preex}
\emph{Let 
\[
g(\z,z) = 1 - \nabla \frac{1}{2\pi i}\frac{\bar \z \cdot d\eta}{1 + |\z |^2}
= \frac{1+ \bar \z \cdot z}{1 + |\z |^2} - \dbar \frac{i}{2\pi}
\frac{\bar \z \cdot d\eta}{1 + |\z |^2}, 
\]
then $g$ is a weight for all $(\z,z)$. Take a $(p,q)$-form
$\phi(\z)$ which grows 
polynomially as $|\z| \to \infty$. If we let $K = (u \w g^k)_{n,n-1}$ and $P
= (g^k)_{n,n}$, then 
\[
\phi(z) =  \int_{|\z| = R} K \w \phi + \int_{|\z| \leq R} K \w \dbar \phi +  
\dbar_z \int_{|\z| \leq R} K \w \phi + \int_{|\z| \leq R} P \w \phi. 
\]
If $k$ is large enough, then the weight will compensate for the growth
of $\phi$, so that the boundary integral will go to zero when $R \to
\infty$. We get the representation 
\[
\phi(z) =  \int K \w \dbar \phi + \dbar_z \int K \w \phi + \int P \w \phi. 
\]}
\end{preex}
Note that if $\phi$ in (\ref{iris}) is a closed form and the 
first and fourth terms of the right hand side 
of Koppelman's formula vanish, we get a 
solution of the $\dbar$-problem for $\phi$. 
Note also that the proof of Koppelman's formula works equally well
over $M \times M$, where $M$ 
is any complex manifold, provided that we can find $K$ and $P$ such
that (\ref{miyazaki}) holds. The purpose of the next section is to
find such $K$ and $P$ in a special type of manifold.

\section{A method for finding weighted Koppelman formulas on
  manifolds} 
\label{yuuri}

We will now describe a method which can be used to find integral 
formulas on manifolds in certain cases, and which is modelled on the
one in the previous section. The method is similar to one presented in
\cite{BE1}, see Remark \ref{kallt} at the end of this section for a
comparison. 

Let $X$ be a complex manifold of dimension $n$, and let $E \to 
X_\z \times X_z$ be a vector bundle of rank $n$, such that we can find
a holomorphic section $\eta$ of $E$ that 
defines the diagonal $\Delta = \{(\z,z): \z=z \}$ of $X \times X$.
In other words, $\eta$ must vanish to the first order on $\Delta$ and be
non-zero elsewhere. Let $\{e_i\}$ be a local frame for $E$, and
$\{e_i^\ast\}$ the dual local frame for $E^\ast$. 
Contraction with $\eta$ is an 
operation on $E^\ast$ which we denote by $\deta$; if $\eta = 
\sum \eta_i e_i$ then
\[
\deta \left(\sum \sigma_i e_i^\ast \right) = \sum \eta_i \sigma_i.
\]
Set 
\[
\nabla_\eta = \deta - \dbar.
\]

Choose a Hermitian metric $h$ for $E$, 
let $D_E$ be the Chern connection on $E$, and $D_{E^\ast}$ the induced
connection on $E^\ast$. Consider $G_E =
  C^\infty (X \times X, \Lambda [T^\ast(X \times X) \oplus E \oplus
  E^\ast])$. If $A$ lies in $C^\infty (X \times X, T^\ast(X \times X)
  \otimes E \otimes E^\ast))$, then we define 
$\tilde{A}$ as the corresponding element in $G_E$, arranged with the
differential form first, then the section of $E$ and finally the
section of $E^\ast$.  For example, if $A =
dz_1 \otimes e_1 \otimes e_1^\ast$, then $\tilde{A} = dz_1 \w e_1
\w e_1^\ast$.  

To define a derivation $D$ on $G_E$, we first let $Df
= \widetilde{D_Ef}$ for a section $f$ of $E$, and $Dg = \widetilde{D_{E^\ast}g}$
for a section $g$ of $E^\ast$. We then extend the definition by 
\[
D(\xi_1 \w \xi_2) = D  \xi_1 \w \xi_2 + (-1)^{\deg\xi_1} \xi_1 \w D \xi_2,
\]
where $D\xi_i = d\xi_i$ if $\xi_i$ happens to be a differential form,
and $\deg\xi_1$ is the total degree of $\xi_1$. For example,
$\deg(\alpha \w e_1 \w e_1^\ast) = \deg\alpha + 2$, where $\deg\alpha$
is the degree of $\alpha$ as a differential form. 
We let 
\[
\mathcal{L}^m = \bigoplus_p C^\infty (X \times X, \Lambda^p E^\ast \w
\Lambda^{p+m} T^\ast_{0,1} (X \times X));
\]
note that $\mathcal{L}^m$ is a subspace of $G_E$. The
operator $\nabla$ will act in a natural way as $\nabla: \mathcal{L}^m
\to \mathcal{L}^{m+1}$. Notice also the analogy 
with the construction (\ref{fisk}) in $\Cn$. As before, if $f \in
\mathcal{L}^m$ and $g \in \mathcal{L}^{k}$, then $f \w g \in
\mathcal{L}^{m+k}$. We also see that $\nabla$ obeys Leibniz' rule, and
that $\nabla^2 = 0$. Let $\End(E)$ denote the bundle of endomorphisms of $E$. 

\begin{prop} \label{flod}
If $v$ is a differential form taking values in $\End(E)$,
and $D_{\End(E)}$ is the induced Chern connection on $\End(E)$, then
\be \label{nana}
\widetilde{D_{\End(E)}v} = D \tilde{v}.
\ee  
\end{prop}

\begin{proof}
Suppose that $v = f \otimes g$, where $f$ is a section of $E$ and $g$
a section of $E^\ast$. We prove first that 
\be \label{pannkaka}
D_{\End(E)}v = D_E f \otimes g + f \otimes D_{E^\ast}g. 
\ee
In fact, if $s$ takes values in $E$, we have 
\[
(D_{\End(E)}v).s = D_E((g.s)f) - (g.(D_E s))f = d(g.s)f +
(g.s)D_Ef - (g.(D_E s))f = 
\]
\[
= (g.s)D_Ef +
(D_{E^\ast}g.s)f = (D_E f \otimes g + f \otimes D_{E^\ast}g).s , 
\]
which proves (\ref{pannkaka}). We have 
\[
\widetilde{D_{\End(E)}v} = \widetilde{D_E f \otimes g} + \widetilde{f
  \otimes D_{E^\ast}g} = Df \w g - f \w Dg = D \tilde{v}
\]
which proves (\ref{nana}). If $v = \alpha \otimes f
\otimes g$, where $\alpha$ is a differential form, we would have
$D_{\End(E)} v = d \alpha \otimes f \otimes g + (-1)^{\deg \alpha}
\alpha \otimes D_{\End(E)} (f \otimes g)$, so the result follows by an
application of $\sim$. Since any differential form taking values in
$\End(E)$ is a sum of such elements, the result follows by linearity. 
\end{proof}

\df{For a form $f(\z,z)$ on $X \times X$, we define
\[
\int_E f(\z,z) \w e_1 \w e_1^\ast \w \ldots \w e_n \w e_n^\ast = f(\z,z).
\]
}
\noindent
Note that if $I$ is the identity on $E$, then $\tilde{I} = e \w
e^\ast = e_1 \w e_1^\ast + \ldots + e_n \w e_n^\ast$. It follows that
$\tilde{I}_n = e_1 \w e_1^\ast \w \ldots \w e_n \w e_n^\ast$ (with the
notation $a_n = a^n/n!$), so the
definition above is independent of the choice of frame. 

\begin{prop} \label{ek}
If $F \in G_E$ then 

\[
d \int_E F = \int_E DF.
\]
\end{prop}

\begin{proof}
If $F = f \w \tilde{I}_n$ we have $d \int_E F = df$ and 
\[
\int_E DF = \int_E [ df \w \tilde{I}_n \pm f \w D(\tilde{I}_n)].
\]
It is obvious that $D_{\End(E)} I = 0$, and by
Proposition \ref{flod} it follows that $D \tilde{I} = 0$, so we are finished.
\end{proof}

We will now construct integral formulas on $X \times X$. As a first
step, we find a section $\sigma$ of $E^\ast$ such that $\deta 
\sigma = 1$ outside $\Delta$. For reasons that will become apparent,
we choose $\sigma$ to have minimal pointwise norm with respect to the metric
$h$, which means that $\sigma = \sum_{ij} h_{ij} \bar 
\eta_j e^\ast_i /|\eta|^2$. Close to $\Delta$, it is obvious that
$|\sigma| \lesssim 1/|\eta|$, and a calculation shows that we also
have $|\dbar \sigma| \lesssim 1/|\eta|^2$. Next, we construct a
section $u$ with the property that $\nabla u = 1 
- R$ where $R$ has support on $\Delta$. We set 

\be \label{nord}
u = \frac{\sigma}{\nabla_\eta \sigma} = \sum_{k=0}^\infty \sigma 
\w (\dbar \sigma)^k,
\ee
note that $u \in \mathcal{L}^{-1}$. By $u_{k,k-1}$ we will mean
the term in $u$ with degree $k$ in $E^\ast$ and degree $k-1$ in
$T^\ast_{0,1} (X \times X)$. It is easily checked that $\nabla
u = 1$ outside $\Delta$.
 
We will need the following lemma:

\begin{lma}\label{dimma}
If $\Theta$ is the Chern curvature tensor of $E$, then 
\[
\nabla_\eta \left(\frac{D \eta}{2 \pi i} + \frac{i \tilde{\Theta}}{2
    \pi}\right) = 0.  
\]
\end{lma}

\begin{proof}
The lemma will follow from the more general statement that if $v$
takes values in 
$\End(E)$, then $\deta \tilde{v} = -v . \eta$. In fact, let $v = f
\otimes g$, where $f$ is a section of $E$ and $g$ 
a section of $E^\ast$, then we have $\deta (f \w g) = -f \w \eta . g =
- (f \otimes g).\eta$. 
Now, note that $\dbar \tilde{\Theta} = 0$ since $D$ is the Chern
connection. We have 
\[
\nabla_\eta \left(\frac{D \eta}{2 \pi i} + \frac{i \tilde{\Theta}}{2
    \pi}\right) =
- \frac{1}{2 \pi i} \left[\dbar D \eta + \deta \tilde{\Theta} \right]= -
  \frac{1}{2 \pi i} \left[\Theta \eta - \Theta \eta
  \right] = 0.  
\]
In the calculations we use that $\eta$ is
holomorphic and that $\dbar \theta = \Theta$ where $\theta$ is the
connection matrix of $D_E$ with respect to the frame $e$.
\end{proof}
\noindent
The following theorem yields a Koppelman formula by Theorem \ref{koala}:
\begin{thm} \label{godis}
Let $E \to X \times X$ be a vector bundle with a section $\eta$ which
defines the diagonal $\Delta$ of $X \times X$. 
We have 

\[
\dbar K = [\Delta] - P,
\]
where 
\be \label{calopteryx}
K = \int_E u \w \left(\frac{D \eta}{2 \pi i} + \frac{i
    \tilde{\Theta}}{2 \pi}\right)_n 
\quad \textrm{and} \quad P = \int_E \left(\frac{D \eta}{2 \pi i} + \frac{i}{2 \pi}\tilde{\Theta}\right)_n, 
\ee
and $u$ is defined by (\ref{nord}). 
\end{thm}
\noindent
Note that since $D \eta$ contains no $e_i$'s, we have 
\[
P = \int_E (\frac{i}{2
\pi}\tilde{\Theta})_n = \det\frac{i\Theta}{2 \pi} = c_n(E),
\]
i.\ e.\ the $n$:th Chern class of $E$. 
\begin{proof} We claim that 

\be \label{kajflock}
\frac{1}{(2 \pi i)^n} \int_E R \w (D \eta)_n = [ \Delta ],
\ee
where $R$ is defined by $\nabla u = 1-R$. If this were true, we would
have by Lemma \ref{dimma} and Proposition \ref{ek}
\begin{eqnarray} 
& & \dbar \int_E u \w \left(\frac{D \eta}{2 \pi i} + \frac{i
    \tilde{\Theta}}{2 \pi}\right)_n =
\int_E \dbar \left[ u \w \left(\frac{D \eta}{2 \pi i} + \frac{i
    \tilde{\Theta}}{2 \pi} \right)_n \right] = \nonumber \\
& = & - \int_E \nabla \left[ u \w \left(\frac{D \eta}{2 \pi i} + \frac{i
    \tilde{\Theta}}{2 \pi} \right)_n \right] = \nonumber \\
& = & - \int_E \left(\frac{D \eta}{2 \pi i} + \frac{i
    \tilde{\Theta}}{2 \pi}\right)_n + \frac{1}{(2 \pi i)^n} \int_E R 
\w (D \eta)_n = [\Delta] - P. \nonumber 
\end{eqnarray}

We want to use Proposition \ref{isis} to prove the claim
(\ref{kajflock}), so we 
need to express the left hand side of (\ref{kajflock}) in local
coordinates. Since $\eta$ defines $\Delta$, we can choose 
$\eta_1, \ldots,
\eta_n$ together with some functions $\tau_1, \ldots, \tau_n$ to 
form a coordinate system locally in a neighborhood of $\Delta$. We have

\[
\frac{1}{(2 \pi i)^n} \int_E R \w (D \eta)_n = \dbar \frac{1}{(2 \pi i)^n}
\int_E \sigma \w (\dbar \sigma)^{n-1} \w (D
\eta)_n, 
\]
and 
\[
\int_E \sigma \w (\dbar \sigma)^{n-1} \w (D\eta)_n = s \w (\dbar s)^{n-1}
+ A,
\]
where $s = \sum \sigma_i d \eta_i$ and $A$ contains only terms which
lack some $d\eta_i$, i.\ e., every
term in $A$ will contain at least one $\eta_i$. Note that both $s$ and $A$ are
now forms in $\Cn$. Recall that we have $|\sigma| \lesssim 1/|\eta|$
and $|\dbar \sigma| 
\lesssim 1/|\eta|^2$ close to $\Delta$ (this is why we chose $\sigma$
to have minimal norm). Thus, by Theorem \ref{isis} we know that 
\[
\dbar [s \w (\dbar s)^{n-1}] = [\Delta],
\]
so it suffices to show that $\dbar A = 0$ in the current sense. But since every
term in $A$ contains at least one $\eta_i$, the singularities which
come from the $\sigma_i$'s and $\dbar \sigma_i$'s will be
alleviated, and in fact we have $A = 
\mathcal{O}(|\eta|^{-2n+2})$. A calculation shows also that $\dbar A = 
\mathcal{O}(|\eta|^{-2n+1})$, and it follows that $\dbar A = 0$ (also cf
the proof of Proposition \ref{isis}). 
\end{proof}
\noindent
It should be obvious from the proof that instead of $u = \sigma/\nabla
\sigma$, we can choose any $u$ such that $\nabla u = 1$ outside
$\Delta$ and $|u_{k,k-1}| \lesssim |\eta|^{-2k+1}$. 

We will obtain more flexible formulas if we use weights:

\df{The section $g$ with values in $\mathcal{L}_0$ is a weight if 
$\nabla g = 0$ and $g_{0,0}(z,z) = 1$.}
\noindent
Theorem \ref{godis} goes through with essentially the same proof if we
take 
\be \label{thai}
K = \int_E u \w g \w \left(\frac{D \eta}{2 \pi i} + \frac{i
    \tilde{\Theta}}{2 \pi}\right)_n
\quad \textrm{and} \quad P = \int_E g \w \left(\frac{D \eta}{2 \pi
    i} + \frac{i \tilde{\Theta}}{2 \pi}\right)_n, 
\ee
as shown by the following calculation: 
\[
\dbar K = -\int_E \nabla u \w g \w \left(\frac{D \eta}{2 \pi i} + \frac{i
    \tilde{\Theta}}{2 \pi}\right)_n = 
-\int_E (g - R) \w \left(\frac{D \eta}{2 \pi i} + \frac{i
    \tilde{\Theta}}{2 \pi}\right)_n = [\Delta ] - P,
\]
which follows from the proof of Theorem \ref{godis} and the properties of
weights. In the next section we will make use of weighted formulas.

If $L$ is a line bundle over $X$, let $L_\z$ denote the line bundle
over $X_\z \times X_z$ defined by $\pi^{-1}(L)$ where $\pi: X_\z
\times X_z \to X_\z$. If we want to find formulas for $(p,q)$-forms
$\phi(\z)$ taking 
values in some line bundle $L$ over $X$, we can use a weight $g$
taking values in $L_z \otimes L_{\z}^\ast$. In fact, then $K$ and $P$
will also take values in $L_z \otimes L_{\z}^\ast$, so that $\phi \w
K$ and $\phi \w P$ take values in $L_z$. Integrating over $\z$, we
obtain $\phi(z)$ taking values in $L$. 

\begin{preremark} \label{kallt}
\emph{To obtain more general formulas, one can find forms $K$ and $P$ such
that 
\be \label{crowley}
dK = [\Delta] - P
\ee
by setting $\nabla'_\eta = \deta - D$ and
checking that the corresponding Lemma \ref{dimma} and Theorem
\ref{godis} are still 
valid. The main difference lies in the fact that since $(\nabla')^2 \neq
0$, we do not have $\nabla' u = 1$ outside $\Delta$, but rather 
\[
\nabla' u = 1 - \frac{\sigma}{(\nabla' \sigma)^2} \w (\nabla')^2 \sigma.
\]
A calculation shows that $(\nabla')^2 \sigma = \delta_\sigma (D \eta -
\tilde{\Theta})$, where $\delta_\sigma$ operates on sections of
$E$. We have $\delta_\sigma (D \eta - \tilde{\Theta}) \w (D \eta -
\tilde{\Theta})^n = \delta_\sigma (D \eta - \tilde{\Theta})^{n+1} = 0$
for degree reasons, so that Theorem \ref{godis} will still hold with
$\nabla$ replaced by $\nabla'$. We can use weights in the same
way, if we require that a weight $g$ has the property $\nabla' g =
0$ instead of $\nabla g = 0$. In this article we are interested
in applications which only require the formulas obtained by using
$\nabla$. }

\emph{In \cite{BE1} Berndtsson obtains $P$
and $K$ satisfying (\ref{crowley}) by a different means, resulting
in the same formulas, but without weights. Also noteworthy is that
$\nabla'$ is a 
superconnection in the sense of Quillen \cite{QU}, and our $\nabla$ is the
$(0,1)$-part of this superconnection. Lemma \ref{dimma} for $\nabla'$
is a Bianchi identity for the superconnection. }
\end{preremark}

\section{Weighted Koppelman formulas on $\Pn$} \label{svamp}

We will now apply the method of the previous section to $X = \Pn$. 
We let $[\z] \in \Pn$ denote the equivalence class of $\z \in
\mathbb{C}^{n+1}$. 
In order to construct the bundle $E$, we first let $F' =  
\mathbb{C}^{n+1} \times (\Pn_{[\z]} \times \Pn_{[z]})$ be the trivial
bundle of rank $n+1$ over $\Pn_{[\z]} \times \Pn_{[z]}$. We next let
$F$ be the bundle of rank $n$ over 
$\Pn_{[\z]} \times \Pn_{[z]}$ which has the fiber
$\mathbb{C}^{n+1}/(\z)$ at the point $([\z],[z])$; $F$ is thus a
quotient bundle of $F'$. If $\alpha$ is a section of $F'$,
we denote its equivalence class in $F$ with $[\alpha]$. We
will not always bother with writing out the brackets, since it will
usually be clear from the context whether a section is to be seen as
taking values in $F'$ or $F$. Let $L^{-1}$ denote
the tautological line bundle of $\Pn$, that is, 
\[
L^{-1} = \{([\z],\xi) \in \Pn \times \mathbb{C}^{n+1}: \xi \in
\mathbb{C} \cdot \z\}
\]
We also define $L^{-k} = (L^{-1})^{\otimes k}$, $L^1 = (L^{-1})^\ast$
and $L^{k} = (L^{1})^{\otimes k}$.
Finally, let $E = F \otimes L^1_{[z]} \to \Pn_{[\z]} \times
\Pn_{[z]}$. Observe that $E$ is thus a subbundle of $E' = F' \otimes
L^1_{[z]}$. It follows that $E^\ast = F^\ast \otimes L^{-1}_{[z]}$,
where $F^\ast = \{ \xi \in (F')^\ast :\xi \cdot \z = 0\}$.  
Berndtsson has the same setup in Example 3, page 337 of
\cite{BE1}, but does not develop it as much (cf Remark \ref{kallt} above). 

Let $\{e_i\}$ be an orthonormal basis of  
$F'$. The section $\eta$ (cf Section \ref{yuuri}) will be $\eta = z
\cdot e = z_0e_0 + \ldots z_ne_n$. Note that $\eta$ takes values  
in $(F') \otimes L^1_{[z]}$, and will thus define an equivalence
class in $F \otimes L^1_{[z]} = E$. The section $\eta$
defines the diagonal since $[\eta(\z,\z)] = [\z \cdot e] = [0]$, so that
$\eta$ vanishes to the first order on $\Delta$.   

We will now choose a metric on $E$. On $F'$ we choose the trivial
metric, which induces the trivial metric also on $(F')^\ast$ and
$F^\ast$. For $[\omega]$ taking values in $F = 
F'/(\z)$, the metric induced from $F'$ is $\| [\omega] \|_F = \| \omega - \pi
\omega\|_{F'}$, where $\pi$ is the orthogonal projection $F' \to
(\z)$.  
Since sections of $L^1_{[z]}$ can be seen as polynomials on
$\mathbb{C}^{n+1}$ which are $1$-homogeneous in $z$, we choose the
metric on $E = F \otimes L^1_{[z]}$ to be 
\be \label{stym}
\| \alpha \otimes [\omega] \|_E = \| \omega - \pi
\omega\|_{F'}|\alpha|/|z|
\ee
for $\alpha \otimes [\omega] \in E$.
We introduce the notation $\alpha \cdot \gamma := \alpha_1 \w
\gamma_1 + \ldots 
\alpha_1 \w \gamma_1$, where $\alpha$ and $\gamma$ are tuples
containing differential forms or sections of a bundle. 

\begin{prop}
Let $\omega \cdot e$ be a section of $E$. The Chern connection and
curvature of $E$ are
\be \label{troll}
D_E (\omega \cdot e) = d \omega \cdot e - \frac{d\z \cdot
  e}{|\z|^2} \w \bar \z \cdot \omega - \partial \log |z|^2 \w \omega \cdot e 
\ee
\be \label{sommar}
\tilde \Theta_E = \ddbar \log |z|^2 \w e^\ast \cdot e - \db
\frac{\bar \z \cdot e^\ast}{|\z|^2} \w d \z \cdot e, 
\ee
with respect to the metric (\ref{stym}) and expressed in the frame $\{e_i\}$
for $F'$.
\end{prop}

\begin{proof}
We begin with finding $D_F$. Let
$\hat \omega \cdot e = (\omega \cdot \bar \z /|\z|^2) \z \cdot e$ be the
projection of $\omega \cdot e$ onto $(\z \cdot e)$. Since the Chern
connection $D_{F'}$ on $F'$ is just $d$, it is easy to see that $D_F
[\omega \cdot e] = [d(\omega \cdot e- \hat \omega \cdot e)]$. We have 
\[
D_F [\omega \cdot e] = [d(\omega \cdot e - \hat \omega \cdot e)] = [d
\omega \cdot e - \frac{d\z \cdot e}{|\z|^2} \w \bar \z \cdot \omega],
\]
since if $d$ does not fall on $\z$ in the second term we get something
that is in the zero equivalence class in $F$. If $\omega \cdot e$ is
projective to start with, 
so will $d\omega \cdot e$ be, and $d\z \cdot e$ is a projective form
since $\delta_{\z} (d\z \cdot e) = \z \cdot e = 0$ in $F$. 

Since the metric on $L^1_{[z]}$ in the local frame $z_0$ is
$|z_0|^2/|z|^2$, the local connection matrix will be $\partial \log
(|z_0|^2/|z|^2)$. If $\xi$ takes values in $L^1_{[z]}$, we get 
\[
D_{L^1_{[z]}} \xi = [d(\xi/z_0) + \partial \log (|z_0|^2/|z|^2)
\xi/z_0] z_0 = d\xi - \partial \log |z|^2 \xi.
\]
It is easy to see that $d(\xi/z_0) + \partial \log (|z_0|^2/|z|^2)
\xi/z_0$ is a projective form, so $d\xi - \partial \log |z|^2 \xi$ is
also projective. Combining the contributions from $L^1_{[z]}$ and $F$,
we get (\ref{troll}), from which also (\ref{sommar}) follows.  
\end{proof}

We want to find the solution $\sigma$ to the equation $\deta
\sigma = 1$, such that $\sigma$ has minimal norm in $E^\ast$. 
It is easy to see that the section $\bar z \cdot e^\ast /|z|^2$ is the
minimal solution to $\deta v = 1$ in the bundle $(E')^\ast = (F')^\ast \otimes
L^{-1}_{[z]}$. The projection of $\bar z \cdot e^\ast /|z|^2$ onto the
subspace $E^\ast$ is   
\[
s = \frac{\bar z \cdot e^\ast}{|z|^2} - \frac{\bar z \cdot \z}
{|\z|^2|z|^2} \bar \z \cdot e^\ast.
\]
Since $\bar z \cdot e^\ast /|z|^2$ is minimal in $(F')^\ast \otimes
L^{-1}_{[z]}$, $s$ must be the minimal solution in $E^\ast$. 
Finally, we normalize to get $\sigma = s /\delta_\eta s$. According to
the method of the previous section, we can then 
set $u = \sigma/\nabla \sigma$ and obtain the forms $P$ and $K$ which
will give us a Koppelman formula (see Theorem \ref{godis}).  

\begin{preremark}
\emph{In local coordinates, for example where $\z_0,z_0 \neq 0$, we have 
\[
|\eta|^2 = \delta_\eta s = \frac{|\z|^2|z|^2 - |\bar z \cdot
  \z|^2}{|\z|^2|z|^2} = \frac{(1+|\z'|^2)(1+|z'|^2) - |1+\bar z' \cdot
  \z'|^2}{(1+|\z'|^2)(1+|z'|^2)},
\]
where $\z' = (\z_1/\z_0, \ldots, \z_n/\z_0)$ and analogously for
$z'$. For the denominator we locally have $(1+|\z'|^2)(1+|z'|^2) \leq
C$ for some constant $C$. As for the numerator, we have 
\begin{eqnarray*}
& & (1+|\z'|^2)(1+|z'|^2) - |1+\bar z' \cdot \z'|^2 =\\
& = & 1 + |\z'|^2 + |z'|^2 + |\z'|^2|z'|^2 - (1 + 2 Re |\bar z' \cdot
\z'| + |\bar z' \cdot \z'|^2) =\\
& = & |z'-\z'|^2 + |\z'|^2|z'|^2 - |\bar z' \cdot \z'|^2 \geq
|z'-\z'|^2
\end{eqnarray*}
In all, we have $\delta_\eta s \gtrsim |z'-\z'|^2$.}
\end{preremark}

To compute integrals of the type (\ref{thai}), we need the 
following proposition.

\begin{lma} \label{kryptogam}
Let $A \stackrel{i}{\hookrightarrow} A'$, where $A'$ is a given 
vector bundle with a given metric and $A = \{\xi \textrm{ taking
  values in } A':  
f \cdot \xi = 0\}$ for a fixed $f$ taking values in $(A')^\ast$. Let
$s$ be the dual section to $f$, and $\pi$ be the orthogonal projection
$\pi: G_{A'} \to G_A$ induced by the metric on $A$. If $B' \in
G_{A'}$, and $B = \pi B'$, then  

\[
\int_A B = \int_{A'} f \w s \w B'.
\] 
\end{lma}

\begin{proof}
We can choose a frame for $A'$ so that $e_0 = s$, and then extend it to an ON
frame for $A'$, so that $A = \Span(e_1, \ldots, e_n)$. 
If we set $e_0^\ast = f$, we have  
\[
\int_{A'} f \w s \w B' = \int_{A'} e_0 \w e_0^\ast \w
\pi B' = \int_A B
\]
and we are done, since the integrals are independent of the frame. 
\end{proof}
\noindent
Note that if $E = A \otimes L$, where $L$ is a line bundle, and $B \in
G_E$, then $\int_E B = \int_A B$. At least, this is true if we
interpret the latter integral to mean that if $g$ is a local frame for
$L$ and $g^\ast$ a local frame for $L^\ast$, then $g$ and $g^\ast$
should cancel out. Since there are as many elements from $L$ as there
are from $L^\ast$, there will be no line bundle elements left. 

We will apply Lemma \ref{kryptogam} with $A = E$, $A' =
E'$ and $f = \z \cdot e^\ast$. We then have
\[
P = \int_{E} \left(\frac{D \eta}{2 \pi i} + \frac{i
    \tilde{\Theta}}{2 \pi}\right)_n = \int_{E'} \frac{\bar \z
  \cdot e \w \z \cdot e^\ast}{|\z|^2} \w \left(\frac{D \eta}{2 \pi i} + \frac{i
    \tilde{\Theta}}{2 \pi}\right)_n 
\]
and similarly for $K$ (this makes it easier to write down $P$ and $K$
explicitly). 

By Theorem \ref{godis}, we have

\[
\dbar K = [\Delta] - P.
\]
(These $K$ and $P$ are also found at the very end of \cite{BE1}.) We
will now modify the method slightly, since in
the paper \cite{LIC} we found formulas for $(0,q)$-forms (derived in a
slightly different way) which are
more appealing than those we have just found, in that we get better
results when we use them to solve 
$\dbar$-equations. We would thus like to have
formulas for $(p,q)$-forms that  
coincide with those of \cite{LIC} in the $(0,q)$-case.

The bundle $F^\ast$ is actually isomorphic to 
$T^\ast_{1,0}(\Pn_{[\z]})$, and an explicit isomorphism is given by $\beta = 
d\z \cdot e$. In fact, if $\xi \cdot e^\ast$ takes values in $F^\ast$, then
$\beta(\xi) = d\z \cdot \xi$. Since $\xi \cdot \z = 0$, the 
contraction of $\beta(\xi)$ with the vector field $\z \cdot 
\partial/\partial \z$ will be zero, so $\beta(\xi) \in
T^\ast_{1,0}(\Pn_{[\z]})$. If $v_{e^\ast}$ is a form with values in
$\Lambda^n E^\ast$, then it is easy to see that 

\be \label{irriterad}
\int_E v_{e^\ast} \w \beta_n = v_{d\z},
\ee
where we get $v_{d\z}$ by replacing every instance of $e_i^\ast$ in
$v_{e^\ast}$ with $d\z_i$. For example, if $v_{e^\ast} = f(\z,z)e^\ast_0 \w \ldots 
\w e^\ast_n$, then $v_{d\z} = f(\z,z)d\z_0 \w \ldots \w d\z_n$. We can
use this to construct integral formulas for  
$(0,q)$-forms with values in $L_{[\z]}^{-n}$, by setting 

\[
K = \int_E u \w \beta_n. 
\]
The formulas we get from this are the same as in \cite{LIC}. We will
now combine these formulas with the ones in (\ref{calopteryx}):

\begin{thm} \label{virvel}
Let $D \subset \Pn$. If $\phi(\z)$ is a $(p,q)$-form 
with values in $L_{[\z]}^{-n+p}$ and 

\begin{eqnarray} \label{arrgh}
K_{p} & = & \int_E u \w \beta_{n-p} \w \left(\frac{D \eta}{2 \pi i} + \frac{i
    \tilde{\Theta}}{2 \pi}\right)_p,\\
P_{p} & = & \int_{e} \beta_{n-p} \w \left(\frac{D \eta}{2 \pi i} + \frac{i
    \tilde{\Theta}}{2 \pi}\right)_p, \nonumber
\end{eqnarray}
with $\beta = d\z \cdot e^\ast$, 
we have the Koppelman formula
\[
\phi([z]) = \int_{\partial D} \phi K_{p} \w \phi +\int_{D} \dbar K_{p}
\w \phi + \dbar_{[z]} \int_{D} K_{p} \w \phi+ \int_{D} P_{p} \w \phi, 
\]
where the integrals are taken over the $[\z]$ variable. 
\end{thm}

\begin{proof} 
We have
\be \label{slash}
\int_E \dbar u \w  \beta_{n-p} \w \left(\frac{D \eta}{2 \pi
    i} + \frac{i \tilde{\Theta}}{2 \pi}\right)_p = [\Delta],
\ee
where $[\Delta]$ should be integrated against sections of $L^{-n+p}$ with
bidegree $(p,q)$. This follows from the proof of Theorem \ref{godis},
since the singularity at $\Delta$ comes only from $u$, and is not
affected by exchanging $\left(\frac{D \eta}{2 \pi i} + \frac{i
    \tilde{\Theta}}{2 \pi}\right)_{n-p}$ for $\beta_{n-p}$.

Using (\ref{slash}), we get
\begin{eqnarray*}
dK_{p} & = & -\int_E \nabla \left[u \w \beta_{n-p} \w \left(\frac{D \eta}{2
    \pi i} + \frac{i \tilde{\Theta}}{2 \pi}\right)_p \right] \\
& = & -\int_E (\nabla u) \w \beta_{n-p} \w \left(\frac{D \eta}{2 \pi
    i} + \frac{i \tilde{\Theta}}{2 \pi}\right)_p = [\Delta] - P_{p}.
\end{eqnarray*}
The Koppelman formula then
follows as in Theorem \ref{koala}.  
\end{proof}

To get formulas for other line bundles, we need to use weights (as
defined in the previous section). We will use the weight 

\[
\alpha = \frac{z \cdot \bar \z}{|\z|^2} - 2 \pi i \dbar \frac{\bar \z 
\cdot e^\ast}{|\z|^2},
\]
note that the first term in $\alpha$ takes values in $L^1_{[z]} \otimes
L^{-1}_{[\z]}$, and the second is a projective form. We then get a
Koppelman formula for $(p,q)$-forms $\phi$ with values in $L^r$ by using 

\begin{eqnarray*}
K_{p,r} & = & \int_E u \w \alpha^{p-n+r} \w \beta_{n-p} \w
\left(\frac{D \eta}{2 \pi 
    i} + \frac{i \tilde{\Theta}}{2 \pi}\right)_p ,\\
P_{p,r} & = & \int_{e}  \alpha^{p-n+r} \w \beta_{n-p} \w \left(\frac{D
    \eta}{2 \pi i} + \frac{i \tilde{\Theta}}{2 \pi}\right)_p . 
\end{eqnarray*}


\begin{preremark} \label{koppel}
\emph{Let $\phi$ be a $(p,q)$-form. Since we cannot raise $\alpha$ to a
negative power, one could wonder how we can get a Koppelman formula if
$\phi$ takes values in $L^r$ where $r < p-n$? In fact, if we look at
the proof of the Koppelman formula in Proposition \ref{koala}, we see
that the roles of $\phi$ and $\psi$ are symmetrical: we could just as 
well use the proof to get a Koppelman formula for the $(n-p,n-q)$-form
$\psi$ which takes values in $L^{-r}$, using the kernels $K_{p,r}$ and
$P_{p,r}$ in Theorem \ref{virvel}. This is a concrete realization of
Serre duality, which in our case says that 
\[
H^{p,q}(\Pn, L^r) \simeq H^{n-p,n-q}(\Pn, L^{-r}).
\]
We will make use of this dual technique when we look at
cohomology groups in the next section.}
\end{preremark}

\begin{preremark}
\emph{In \cite{BE2} Berndtsson constructs integral formulas for sections of
line bundles over $\Pn$. These formulas coincide with ours in the case
$p=0$, but they are not the same in the general
$(p,q)$-case. Nonetheless, they do give the same result as our
formulas when used to find the trivial cohomology groups of the line
bundles of $\Pn$ (see the next section). More precisely, his 
formulas can also be used to prove Proposition \ref{lunch} below, but
no more, at least not in any obvious way. } 
\end{preremark}

\section{An application: the cohomology of the line bundles of $\Pn$}
\label{wombat}

Let $D$ in Theorem \ref{virvel} be the whole of $\Pn$; then
the boundary integral will disappear. The only obstruction to solving
the $\dbar$-equation is then the term containing $P_{p,r}$.
We will use our explicit formula for $P_{p,r}$ 
to look at the cohomology 
groups of $(p,q)$-forms with values in different line bundles, 
and determine which of them are trivial. We have 
\begin{eqnarray*} 
& P_{p,r} = & \int_{e} \beta_{n-p} \w \left(\frac{D \eta}{2 \pi
    i} + \frac{i \tilde{\Theta}}{2 \pi}\right)_p \w 
\alpha^{n-p+r} = \\
& = & \int_{e'} \frac{\bar \z \cdot e \w \z \cdot e^\ast}{|\z|^2} \w 
(d\z \cdot e)_{n-p} \w 
\left(dz \cdot e - \frac{z \cdot \bar \z}{|\z|^2}\w d\z \cdot 
e -  \right. \\
& - & \left. \frac{\partial |z|^2}{|z|^2} z \cdot e +\omega_z e^\ast 
\cdot e - \frac{d \bar 
\z \cdot e^\ast \w d \z \cdot e}{|\z|^2}\right)_p \w \left( \frac{z 
\cdot \bar \z}
{|\z|^2} - \dbar \frac{\bar \z \cdot e^\ast}{|\z|^2} \right)^{n-p+r}.
\end{eqnarray*}
We can now prove:

\begin{prop} \label{lunch}
From the formula for $P_{p,r}$ just above, it follows that the cohomology groups
$H^{p,q}(\Pn, L^r)$ are trivial in the following cases:\\
\\
\noindent
a) $q = p \neq 0,n$ and $r \neq 0$.\\
b) $q = 0$, $r \leq p$ and $(r,p) \neq (0,0)$.\\
c) $q = n$, $r \geq p-n$ and $(r,p) \neq (0,n)$. \\
d) $p < q$ and $r \geq -(n-p)$.\\
e) $p > q$ and $r \leq p$.\\
\end{prop}
\noindent
Unfortunately, these are not all the trivial cohomology groups;
instead of d) and e) we should ideally get that the groups are trivial
for $q \neq 0,n,p$ (cf \cite{DE} page 397). 

\begin{proof} The general strategy is this: we take a 
$\dbar$-closed form $\phi(z)$ of given 
bidegree and with values in a given line bundle, and then try to show
that $\phi(z)$ is  
exact by means of the Koppelman formula. One possibility of 
doing this is proving that $\int_\z 
\phi(\z) \w P_{p,r}(\z,z) = 0$, which can be either because the 
integrand is zero, or because the 
integrand is $\dbar_\z$-exact (since then Stokes' formula can 
be applied). 
Another possibility is proving that $P_{p,r}$ is $\dbar_z$-exact, 
since then $\int_\z \phi \w P_{p,r}$ 
will be $\dbar_z$-exact as well. \\
\\
\noindent
\textbf{Proof of a):} Let $r > 0$ and $p = q \neq 0,n$; we must then
look at the term in $P_{p,r}$ with bidegree $(p,p)$ in 
$z$ and $(n-p,n-p)$ in $\z$, it is equal to  

\be \label{traktor}
\int_{e'} \frac{\bar \z \cdot e \w \z \cdot e^\ast}
{|\z|^2} \w (d\z \cdot e)_{n-p} \w 
\left( \omega_z \w e^\ast \cdot e \right)^p \w
\left( \frac{z \cdot \bar \z} {|\z|^2}\right)^{r} \w
\left( \frac{d \bar \z \cdot e^\ast}{|\z|^2} \right)^{n-p}  \! \! \!
\! \! \! \! \! \!. 
\ee
We will show that (\ref{traktor}) is actually $\dbar_z$-exact.
The factor in (\ref{traktor}) which depends on $z$ is $(z \cdot 
\bar \z)^{r} \omega_z^p$, which is at 
least a $\dbar_z$-closed form. Can we write $(z \cdot \bar \z)^{r} 
\omega_z^p = \dbar_z g(z)$, where
$g$ is a projective form? Actually, we have $\dbar_z [(\bar \z 
\cdot z)^r \partial |z|^2/|z|^2 \w 
\omega_z^{p-1}] = (z \cdot \bar \z)^r \omega_z^p$, but $(\bar \z 
\cdot z)^r \partial |z|^2/|z|^2 \w 
\omega_z^{p-1}$ is not a projective form. This can be remedied by 
adding a holomorphic term 
$(\bar \z \cdot z)^{r-1} (\bar \z \cdot dz) \w \omega_z^{p-1}$, 
since then we can take 

\[
g = (\bar \z \cdot z)^{r-1} [(\bar \z \cdot z) \frac{\partial |z|^2}
{|z|^2}  - \bar \z 
\cdot dz] \w \omega_z^{p-1}.
\]
Since (\ref{traktor}) is $\dbar_z$-exact, we have proved a) when $r>0$. 
If $-r<0$, by Remark \ref{koppel} in the previous section we must look
at $P_{n-p,r}$,  
which is again $\dbar_z$-exact, and then $\int_z \phi(z) \w 
P_{n-p,r} = 0$ by Stokes' Theorem. \\ 
\\
\noindent
\textbf{Proof of b):} Note that here we really want to prove 
that $\phi = 0$, since $\phi$ cannot 
be $\dbar$-exact. To prove this we again use the dual case in Remark
\ref{koppel}. 
We want to show that $\int_z \phi(z) \w 
P_{n-p,r}(\z,z) = 0$, when $\phi(z)$ has bidegree $(p,0)$ and 
takes values in $L^{-r}_z$. First assume that $p>0$, 
then we must look at the term in $P_{n-p,r}$ of bidegree $(n-p,n)$ in
$z$. No term in $P_{n-p,r}$ has a 
higher degree in $d \bar z$ than in $dz$, so $\int \phi(z) \w 
P_{n-p,r}(\z,z) = 0$. If $p=0$, then then we must look at the term in
$P_{n,r}$ with bidegree $(n,n)$ in $z$ and $(0,0)$ in $\z$. 
The $z$-dependent factor of this term is $(z \cdot \bar \z)^r
\omega_z^n$, which is $\dbar_z$-exact in the same way as in 
the proof of a). This proves the case $p=0$, $-r <0$, but the 
proof breaks down when $r=0$, where there is a non-trivial cohomology. \\
\\
\noindent
\textbf{Proof of c):} First let $p < n$. There is no term in $P_{p,r}$ 
with bidegree $(p,n)$ in $z$, since there are not enough $d \bar z$'s,
so $\int_\z \phi(\z) \w P_{p,r}(\z,z) = 0$.
If $p=n$, we look at the term in $P_{p,r}$ with bidegree $(n,n)$ in
$z$ and $(0,0)$ in $\z$. This is dealt with  
exactly as the case $p=0$ in the proof of b).\\
\\
\noindent
\textbf{Proof of d) and e):} Let $q \neq 0,n,p$. If $p < q$ and $r \geq
-(n-p)$, we look at the term in $P_r$ with bidegree $(p,q)$ in $z$. 
It is zero, since we cannot have more $d \bar z$'s 
than $dz$'s, so $\int_\z \phi(\z) \w P_{p,r} = 0$. 
Similarly, if $p > q$ we use the dual method: the term in 
$P_{n-p,r}$ with bidegree $(n-p,n-q)$ in $z$ is zero when $n-p < n-q$
and $r \geq -p$, again since we cannot have more $d \bar z$'s 
than $dz$'s. This shows that $\int_z \phi(z) \w P_{n-p,r} = 0$ for 
$r \geq -p$, where $\phi$ takes values in $L^{-r}$ and $-r \leq p$. 
\end{proof}

\section{Weighted Koppelman formulas on $\Pn \times \Pm$}
\label{sol}

We will now find integral formulas on $\Pn \times \Pm$. Let 
$([\z], [\tilde{\z}], [z],[\tilde{z}])$ be a point in $(\Pn \times
\Pm) \times (\Pn \times \Pm)$. The procedure
will be quite similar to that of Section \ref{svamp}, but for
simplicity we will limit ourselves to the case of $(0,q)$-forms. This
corresponds to using only $\beta$ in the formula
(\ref{arrgh}). According to formula (\ref{irriterad}), then, we
can construct our kernel directly, without any need to refer to the 
bundle $E$, in the following way (also see \cite{LIC}). 
Let $\eta_\z = 2\pi i z \cdot 
\frac{\partial}{\partial \z}$ and $\eta = \eta_\z +
\eta_{\tilde{\z}}$.  
We take $\deta$ to be contraction with $\eta$ and set 
$\nabla = \deta -\dbar$. Note that $\eta = 0$ on $\Delta$. Now set 

\[
s_{\z} = \frac{\bar z \cdot d \z}{|z|^2} - \frac{\bar z 
\cdot \z}{|z|^2 |\z|^2} \bar \z \cdot d \z
\]
and then $s = s_{\z} + s_{\tilde{\z}}$. Observe that 
$\deta s$ is a scalar, which is zero only on $\Delta$. 

\begin{prop}
If $u = s / \nabla s$, then $u$ satisfies $\nabla u . 
\phi = (1 - [\Delta]) . \phi$, where $\phi$ is a form of 
bidegree $(n+m,n+m)$ which 
takes values in $L^{-n}_{[\z]} \otimes L^{-m}_{[\tilde{\z}]} 
\otimes L^{n}_{[z]} \otimes L^{m}_{[\tilde{z}]}$ and contains no 
$d\z_i$'s or $d \tilde{\z}_i$'s.
\end{prop}

\begin{proof}
The restriction on $\phi$ is another way of saying that our formulas
only will work for $(0,q)$-forms. The proposition will follow from
Theorem \ref{virvel} if we integrate in $\Pn_{[\z]} \times \Pn_{[z]}$ and
$\Pm_{[\tilde \z]} \times \Pm_{[\tilde z]}$ separately. 
\end{proof}


To obtain weighted formulas, let 

\[
\alpha = \frac{z \cdot \bar \z}{|\z|^2} + 2 \pi i \ddbar 
\log |\z|^2,
\]
and let $\tilde \alpha$ be the corresponding form in $([\tilde{\z}],
[\tilde{z}])$. We have $\nabla \alpha = \nabla \tilde \alpha = 0$, so 

\[
\nabla (\alpha^{n+k} \w \tilde \alpha^{m+l} \w u) = 
\alpha^{n+k} \w \tilde \alpha^{m+l} \w \nabla u = \alpha^{n+k} \w
\tilde \alpha^{m+l} - [\Delta],
\]
where $[\Delta]$ must be integrated against sections 
of $L^{k}_{[\z]} \otimes L^{l}_{[\tilde{z}]}$. The following theorem
follows from Theorem \ref{koala}. 

\begin{thm}
If $K = \alpha^{n+k} \w \tilde \alpha^{m+l} \w u$ and 
$P = \alpha^{n+k} \w \tilde \alpha^{m+l}$ we get the Koppelman formula
\begin{eqnarray*}
\phi([z],[\tilde{z}]) & = & \int_{\partial D} \phi([\z],
[\tilde{\z}]) \w K + \int_{D} \dbar 
\phi([\z],[\tilde{\z}]) \w K + \\
& + & \dbar_z \int_{D} \phi([\z],[\tilde{\z}]) \w K + 
\int_{ D} 
\phi([\z],[\tilde{\z}]) \w P \\
\end{eqnarray*}
for differential forms $\phi([\z], [\tilde{\z}])$ on $\Pn 
\times \Pm$ with bidegree $(0,q)$ which take values in 
$L^{k}_{[\z]} \otimes L^{l}_{[\tilde{\z}]}$. 
\end{thm}

Now assume that $\dbar \phi = 0$. For which $q$, 
$k$ and $l$ is $\phi$ $\dbar$-exact? To 
show that a particular $\phi$ is $\dbar$-exact, we need to 
show that the term $\int_{\Pn \times \Pm} 
\phi([\z]) \w P$ either is zero, or is $\dbar$-exact. Since 
$P$ consists of two factors where one 
depends only on $\z$ and the other only on $\tilde{\z}$, we 
can write 

\be \label{mnium}
\int_{\Pn \times \Pm} \phi([\z],[\tilde{\z}]) \w P = 
\int_{\Pm} \left(\int_{\Pn} \phi([\z],
[\tilde{\z}]) \w \alpha^{n+k} \right) \w \tilde \alpha^{m+l}.
\ee
We get the following theorem:

\begin{prop}
We look at differential forms $\phi([\z], [\tilde{\z}])$ on 
$\Pn_{[\z]} \times \Pm_{[\tilde{\z}]}$ with 
bidegree $(0,q)$, which take values in the line 
bundle $L^{k}_{[\z]} \otimes L^{l}_{[\tilde{\z}]}$. The cohomology 
groups $H^{(0,q)} (\Pn \times \Pm, L^{k}_{[\z]} \otimes L^{l}_
{[\tilde{\z}]})$ are trivial in the following cases:\\
\\
a) $q \neq 0,n,m,n+m$ \\
b) $q = 0$ and $k<0$ or $l<0$\\
c) $q = n$ and $l < 0$ or $k \geq -n$\\
d) $q = m$ and $k < 0$ or $l \geq -m$\\
e) $q = n+m$ and $k \geq -n$ or $l \geq -m$.\\
\end{prop}

\begin{proof}
To determine when (\ref{mnium}) is zero, we use Theorem \ref{lunch}. 
Assume that the form $\phi$ has bidegree $(0,q_1)$ in $\z$ and 
$(0,q_2)$ in $\tilde{\z}$ and $q_1 + q_2 = q$. If, for some $q_1$ and $k$, 
we know that $H^{(0,q_1)} (\Pn, L^{k})$ 
is trivial, this means either that $\int_{[\z]} \phi([\z],[\tilde{\z}]) \w 
P([\z],[z]) = 0$ or that $\int_{[\z]} 
\phi([\z],[\tilde{\z}]) \w P([\z],[z]) = \dbar_z a([z], [\tilde{\z}])$
for some $a([z], [\tilde{\z}])$. In the first case, it follows that the 
expression in (\ref{mnium}) is also zero. In the second case, we get 
\begin{eqnarray*}
& & \int_{\Pm} \left(\int_{\Pn} \phi([\z],[\tilde{\z}]) \w 
\alpha^{n+k} \right) \w \tilde 
\alpha^{m+l} = \dbar_z \int_{\Pm} a([z], [\tilde{\z}]) \w 
\tilde \alpha^{m+l} = \\
& =  & \dbar \int_{\Pm} a([z], [\tilde{\z}]) \w \tilde \alpha^{m+l}
\end{eqnarray*}
since the integrand is holomorphic in $[\tilde z]$. The same holds if
$H^{(0,q_2)} (\Pm, L^{l})$ is trivial.  
The conclusion is that $H^{(0,q_1+q_2)} (\Pn \times \Pm,
L^{k}_{[\z]} \otimes L^{l}_ {[\tilde{\z}]}) = 0$ either when 
$q_1$ and $k$ are such that $H^{(0,q_1)} (\Pn, L^{k}) = 0$, or when
$q_2$ and $l$ are such that $H^{(0,q_2)} (\Pm, L^{l}) = 0$. 

Now, we really have a sum 
\[
\phi = \sum_{q_1 + q_2 = q} \phi_{q_1,q_2}
\]
of terms of the type above. For the cohomology group to be trivial,
we must have $\int \phi_{q_1,q_2} \w P = 0$ for all of them. We know
that $q_2 = q-q_1$. If we have either $0 < q_1 < n$ or $0 < q_2 < m$ then 
$\int \phi_{q_1,q_2} \w P = 0$ according to Theorem \ref{lunch}. The
only ways to avoid this are if $q = q_1 = q_2 = 0$; if $q = q_1 = n$
and $q_2 = 0$; if $q_1 = 0$ and $q = q_2 = m$ or if $q = n+m$ and $q_1
= n$, $q_2 = m$. Then a) - e) follow from Theorem \ref{lunch}. 
\end{proof}

\section{Weighted integral formulas on Stein manifolds}
\label{sjunga}

If $X$ is a Stein manifold it is, in general, impossible to find $E
\to X \times X$ and $\eta$ with the desired properties as described in
Section \ref{yuuri}. What is possible is to find a section $\eta$ of
a bundle $E$ such that $\eta$ has good properties close to $\Delta$,
but then $\eta$ will in general have other zeroes as well. It turns
out that it is possible to work around this and still construct
weighted integral formulas. This section relies on the article
\cite{HL} by Henkin and Leiterer, where such an $\eta$ was constructed. 

More precisely, let $\pi$ be the projection from $X_{\z} \times X_z$
to $X_\z$, and $E = \pi^\ast(T_{1,0}(X_\z))$. 
By Section 2.1 in \cite{HL} we have the result

\begin{thm}
There exists a holomorphic section $\eta$ of $E$ such that $\{\eta = 0
\} = \Delta \cup F$, where $F$ is closed and $\Delta \cap F =
\varnothing$. Close to $\Delta$ we have 
\be \label{fram}
\eta(\z,z) = \sum_1^n [\z_i - z_i +
\mathcal{O}(|\z-z|^2)]e_i. 
\ee
Moreover, there exists a holomorphic function $\phi$ such that
$\phi(z,z)=1$ and $|\phi| \lesssim |\eta|$ on a neighborhood of $F$.  
\end{thm}

We define $\deta$, $\nabla$ etc in the same way as in Section
\ref{regn}. Let $s \in E^\ast$ be the section satisfying $\deta s = 1$
outside $\Delta \cup F$ which has pointwise minimal norm, and 
define $u = s/\nabla s$. If we define
\[
K = \int_E \phi^M u \w \left(\frac{D \eta}{2 \pi i} + \frac{i
    \tilde{\Theta}}{2 \pi}\right)_n
\quad \textrm{and} \quad P = \int_E \phi^M \left(\frac{D \eta}{2 \pi
    i} + \frac{i \tilde{\Theta}}{2 \pi}\right)_n, 
\]
where $M$ is large enough that $\phi^M u$ has no singularities on
$F$, then Theorem \ref{godis} applies and we have $\dbar K = [\Delta] - P$. 
In this way, we recover the formula found in Example 2 of \cite{BE1},
except that our approach also allows for weights. We define weights
in the same way as before (note that $\phi$ is in fact a weight). If
$g$ is a weight, we will get a Koppelman formula with
\be \label{attestera}
K = \int_E \phi^M g \w u \w \left(\frac{D \eta}{2 \pi i} + \frac{i
    \tilde{\Theta}}{2 \pi}\right)_n
\  \textrm{and} \ \ P = \int_E \phi^M g \w \left(\frac{D \eta}{2 \pi
    i} + \frac{i \tilde{\Theta}}{2 \pi}\right)_n \! \!. 
\ee

Note that since
$E$ is a pullback of a bundle on $X_\z$, the connection and curvature
forms of $E$ depend only on $\z$. Hence $P =  c_n(E)$ is bidegree
$(n,n)$ in $\z$, and we have $\int_\z P(\z,z) \w \phi(\z) = 0$
except in the case where $\phi$ has bidegree $(0,0)$. The last term in
the Koppelman formula thus presents no obstruction to solving the
$\dbar$-equation on $X$.  

\begin{preex}
\emph{In \cite{FO} there is an example of weighted formulas on Stein 
manifolds, which we can reformulate to fit into the present
formalism. Let $G \subset X$ be a strictly pseudoconvex domain. By 
Theorem 9 in \cite{FO} we can find a function $\psi$ defined on a
neighborhood $U$ of $G$ which embeds $G$ in a strictly convex set $C
\subset \Cn$. If $\sigma$ is the 
defining function for $C$, then $\rho = \sigma \circ \psi$ is a
strictly plurisubharmonic defining function for $G$. On $U$ we
introduce the weight}  
\[
g(\z,z) = \left(1 - \nabla \frac{\frac{\partial \rho(\z)}{\partial \z}
    \cdot e^\ast}{2 \pi i \rho(\z)}
\right)^{-\alpha} = \left( -\frac{v}{\rho} - \omega
\right)^{-\alpha} 
\]
\emph{where}
\[
v = \frac{\partial \rho(\z)}{\partial \z} \cdot \eta -
  \rho(\z)
\quad \textrm{\emph{and}} \quad 
\omega = \dbar \left[ \frac{\frac{\partial \rho(\z)}{\partial \z}
    \cdot e^\ast}{2 \pi i \rho(\z)} \right]
\]
\emph{Note that $g$ is holomorphic in $z$.  If
$\Re \alpha$ is large enough, then $g(\cdot,\z)$ will be zero on $\partial G$,
since $\sigma(\partial C) = 0$. This implies that if $f$ is a
holomorphic function and $P$ is defined by (\ref{attestera}), we will have}

\[
f(z) = \int_G f(\z) P,
\]
\emph{for $z \in G$, by Koppelman's formula. We also have the estimate}
\[
-\rho(\z) - \rho(z) + \epsilon |\z-z|^2 \leq 2 \Re v(\z,z) \leq 
-\rho(\z) - \rho(z) + c|\z-z|^2,
\]
\emph{where $\epsilon$ and $c$ are positive and real. By means of
  this, we can get results in strictly pseudoconvex domains $G$ in
  Stein manifolds similar to ones which are known in strictly
  pseudoconvex domains in $\Cn$. For example, one can obtain a direct
  proof of the Henkin-Skoda theorem which gives
  $L^1$-estimates on $\partial G$ for solutions of the $\dbar$-equation. }



\end{preex}

\begin{preex}
\emph{We can also solve division problems on $X$. Let $D \subseteq X$ be a
domain, and take $f(\z) = (f_1(\z), \ldots, f_m(\z))$ where $f_i \in
\mathcal{O}(\overline{D})$. Assume that $f$ has no common zeroes in
$D$. We want to solve the division problem $\psi = f \cdot p$ in $D$, where
$\psi$ is a given holomorphic function, by means of integral formulas. }

\emph{By Cartan's Theorem B, we can find $h(\z,z) = (h_1(\z,z), \ldots, 
h_m(\z,z))$, where $h_i$ is a holomorphic section of $E^\ast$, such that $\deta
h_i(\z,z) = \phi(\z,z) (f_i(\z) - f_i(z))$. We set 
\[
g_1(\z,z) = (\phi - \nabla (h \cdot
\sigma(\z))^\mu = (\phi f(z) \cdot \sigma + h \cdot \dbar
\sigma)^\mu,
\]
where $\sigma = \bar f /|f|^2$ and $\mu = \min(m,n+1)$, then $g_1$ is
a weight. Now, $f(z)$ is a factor in $g_1$, since $(h \cdot \dbar
\sigma)^\mu = 0$. In fact, we have $(h \cdot \dbar \sigma)^{n+1} = 0$ for
degree reasons, and $(h \cdot \dbar \sigma)^m = 0$ since $f \cdot
\sigma = 1$ implies $f \cdot \dbar \sigma = 0$, so that $\dbar
\sigma_1, \ldots, \dbar \sigma_m$ are linearly dependent. }

\emph{By the Koppelman formula we have
\[
\psi(z) = \int_{\partial D} \psi \phi^M K + \int_D \psi \phi^M P
\]
where $K$ and $P$ are defined by (\ref{attestera}) using the weight
$g_1$. Since $f(z)$ is a factor in $g_1$, we have $\psi(z) = f(z)
\cdot p(z)$, where $p(z)$ will be holomorphic if $D$ is such that we
can find $u$ holomorphic in $z$ (for example if $D$ is pseudoconvex).}

\end{preex}

\end{document}